\documentclass[11pt]{article}
\usepackage{amssymb, amsfonts, amsmath }
\usepackage{latexsym}

\providecommand{\MR}{\relax\ifhmode\unskip\space\fi MR }

\newtheorem{thm}{Theorem}[section]

\newtheorem{lem}[thm]{Lemma}

\newtheorem{dfn}[thm]{Definition}

\newtheorem{qs}[thm]{Question}

\newcommand{\demo}{ \noindent {\it   Proof. }}
\newcommand{\qed}{$\Box$
}
\newcommand{\aut}{\operatorname{Aut}}

\title{A recursive presentation for Mihailova's subgroup}

\author{O.\ Bogopolski \\ \small{Institute of Mathematics of} \\ \small{Siberian Branch of Russian Academy
of Sciences,} \\ {\small Novosibirsk, Russia} \\ \small{e-mail: groups@math.nsc.ru}\\ \\ E. Ventura \\
\small{Dept. Mat. Apl. III,} \\ \small{Universitat Polit\`{e}cnica de Catalunya,} \\ \small{Manresa, Barcelona, Catalunya}
\\ \small{e-mail: enric.ventura@upc.edu} }

\begin{document}
\maketitle

\begin{abstract}
We give an explicit recursive presentation for Mihailova's subgroup $M(H)$ of $F_n\times F_n$ corresponding to a
finite, concise and Peiffer aspherical presentation $H=\langle x_1,\dots ,x_n\,|\, R_1,\dots ,R_m\rangle$. This
partially answers a question of R.I.~Grigorchuk, \cite[Problem 4.14]{Gr}. As a corollary, we construct a finitely
generated recursively presented orbit undecidable subgroup of ${\text{\rm Aut}}(F_3)$.
\end{abstract}

\section{Introduction}\label{intro}

For all the paper, let $n\geqslant 2$, let $F_n$ be the free group with basis $\{ x_1,\ldots ,x_n \}$, and let
$H=\langle x_1,\ldots ,x_n\,|\, R_1,\ldots ,R_m\rangle$ be a finite presentation of a quotient $H$ of $F_n$ (although
most of what follows will depend on the specific presentation, we shall make the usual abuse of notation which consists
on denoting by $H$ both the group and its given presentation).

K.A.~Mihailova, in her influential paper~\cite{Mih}, associated to the presentation $H$
the \emph{Mihailova subgroup} of $F_n \times F_n$, namely
 $$
M(H) =\{ (w_1 ,w_2)\in F_n\times F_n \mid w_1 =_H w_2 \} \leqslant F_n\times F_n,
 $$
i.e. the subgroup of pairs of words in $F_n$ determining the same element in $H$. It is clear that $(x_i,x_i)$ and
$(1,R_j)$ belong to $M(H)$ for all $i=1,\ldots ,n$ and $j=1,\ldots ,m$, and it is not difficult to see that, in fact,
these pairs generate $M(H)$. The important observation made in~\cite{Mih} says that the membership problem for $M(H)$
in $F_n\times F_n$ is solvable (i.e. there exists an algorithm to decide whether a given $(w_1, w_2)\in F_n\times F_n$
belongs to $M(H)$ or not) if and only if the word problem for $H$ is solvable.

By a result of P.S.~Novikov~\cite{N} and W.W.~Boone~\cite{Boo} (see also~\cite{Bor}), there exist finitely presented
groups with unsolvable word problem. Thus, there also exist finitely generated subgroups of $F_n\times F_n$ with
unsolvable membership problem.

Clearly, $M(H)$ has solvable word problem for every $H$ (because $F_n\times F_n$ also does). In particular, $M(H)$ is
recursively presented. More interestingly, F.J.~Grunewald proved, in~\cite[Theorem B]{Gru}, that if $H$ is infinite
then $M(H)$ cannot be finitely presented. In~\cite{BR}, G.~Baumslag and J.E.~Roseblade completely described the
structure of finitely presented subgroups of $F_n\times F_n$, a result that was later reproved by H.~Short~\cite{Short}
and M.R.~Bridson and D.T.~Wise~\cite{BW}, and that implies Grunewald's result.

In this context, a natural problem is to look for recursive presentations for Mihailova's group $M(H)$, in terms of the
original presentation $H$. This was recently posted as Problem~4.14 in~\cite{Gr} by R.I.~Grigorchuk: ``What kind of
presentations can be obtained for Mihailova's subgroups of $F_n\times F_n$ determined by finite automata?"

The main result in the present paper (Theorem~\ref{pres} below) gives a partial answer to this problem: under
certain technical conditions on the initial $H$ we give an explicit recursive presentation for $M(H)$ with
finitely many generators and a one-parametric family of relations.

\begin{thm}\label{pres}
Let $F_n$ be the free group on $x_1,\dots ,x_n$, and let $H=\langle x_1,\dots ,x_n\,|\, R_1,\dots ,R_m\rangle$ be a
finite, concise and Peiffer aspherical presentation. Then Mihailova's group $M(H)\leqslant F_n\times F_n$ admits the
following presentation
 $$
\big\langle d_1,\dots ,d_n, t_1,\dots ,t_m\, \mid \, [t_j, d^{-1}t_i^{-1}r_i \,d],\,\, [t_i,\,\, \text{\rm root}(r_i
)]\,\,\, (1\leqslant i,j\leqslant m,\,\, d\in D_n)\big\rangle,
 $$
where $D_n$ is the free group with basis $d_1,\dots ,d_n$, where $r_i$ denotes the word in $D_n$ obtained from $R_i$ by
replacing each $x_k$ to $d_k$, and where {\rm root}$(r_i)$ denotes the unique element $s_i \in D_n$ such that $r_i$ is
a positive power of $s_i$ but $s_i$ itself is not a proper power.

In this presentation the elements $d_i$ and $t_j$ correspond, respectively, to the elements $(x_i,x_i)$ and $(1,R_j)$
of $M(H)$.
\end{thm}

As a corollary we deduce the existence of a finitely generated, orbit undecidable subgroup of ${\text{\rm Aut}}(F_3)$
(see~\cite{BMV} for details), which has the recursive presentation given in Theorem~1.1.

The structure of the paper is the following. In Section~2 we recall some definitions and discuss some properties of
concise and Peiffer aspherical presentations that will be used later. In Section~3 we prove Theorem~\ref{pres}. And in
Section~4 we recall the relationship between Mihailova's subgroup and orbit undecidability, recently discovered
in~\cite{BMV}, and deduce the announced corollary (Theorem~\ref{corol}).

\section{Asphericity}\label{asp}

As stated, let $\langle x_1,\ldots ,x_n\,|\, R_1,\ldots ,R_m\rangle$ be a finite presentation. Formally, $R_1,\ldots
,R_m $ is a list of words in the alphabet $\{ x_1,\ldots ,x_n \}^{\pm 1}$ which may contain the trivial element,
possible repetitions, and even possible members conjugated to each other or to the inverse of each other.

A presentation $\langle x_1,\ldots ,x_n\,|\, R_1,\ldots ,R_m\rangle$ is called \emph{concise} if every relation $R_i$
is non-trivial and reduced, and every two relations $R_i$, $R_j$, $i\neq j$, are not conjugate to each other, or to the
inverse of each other. Given an arbitrary finite presentation, $\langle x_1,\ldots ,x_n\,|\, R_1,\ldots ,R_m\rangle$,
one can always reduce the relations and eliminate some of them, to obtain another presentation of the same group, which
is concise. We call this a \emph{concise refinement} of $\langle x_1,\ldots ,x_n\,|\, R_1,\ldots ,R_m\rangle$.

Now, we recall the definition of Peiffer transformations. Consider some elements $U_1,\ldots, U_l\in F_n$, some
relators $R_{i_1},\ldots ,R_{i_l}\in \{ R_1,\ldots ,R_n \}$, and some numbers $\varepsilon_1,\ldots ,\varepsilon_l \in
\{-1,1\}$, such that the equation
 $$
(U_1 R_{i_1}^{\varepsilon_1}U_1^{-1}) \cdots (U_l R_{i_l}^{\varepsilon_l }U_l^{-1})=1
 $$
holds in $F_n$. In this situation, the sequence of elements $(U_1 R_{i_1}^{\varepsilon_1 }U_1^{-1},\ldots ,U_l
R_{i_l}^{\varepsilon_l }U_l^{-1})$ of $F_n$ is called an \emph{identity among relations} of length $l$. For $l=0$ we
have the \emph{empty} identity among relations, $(\, )$.

In such a sequence, let us replace two consecutive terms, say $U_p R_{i_p}^{\varepsilon_p} U_p^{-1}$ and
$U_{p+1}R_{i_{p+1}}^{\varepsilon_{p+1}}U_{p+1}^{-1}$ for some $1\leqslant p\leqslant {l-1}$, by the new ones
$U_{p+1}R_{i_{p+1}}^{\varepsilon_{p+1}}U_{p+1}^{-1}$ and $(U_{p+1}R_{i_{p+1}}^{-\varepsilon_{p+1}}U_{p+1}^{-1}U_p)
R_{i_p}^{\varepsilon_p} (U_p^{-1}U_{p+1}R_{i_{p+1}}^{\varepsilon_{p+1}}U_{p+1}^{-1})$. Since the product of the two old
terms do coincide with that of the two new ones, the new sequence is again an identity among relations. This
transformation is called a \emph{Peiffer transformation of the first kind} or, shortly, an \emph{exchange}.

Suppose now that in the sequence $(U_1 R_{i_1}^{\varepsilon_1}U_1^{-1},\ldots ,U_l R_{i_l}^{\varepsilon_l}U_l^{-1})$
there are two consecutive terms, say $U_p R_{i_p}^{\varepsilon_p} U_p^{-1}$ and
$U_{p+1}R_{i_{p+1}}^{\varepsilon_{p+1}}U_{p+1}^{-1}$ for some $1\leqslant p\leqslant {l-1}$, whose product equals 1.
Then, we can obtain a new identity among relations by just deleting these two terms. This transformation and the
inverse one are called \emph{Peiffer transformations of the second kind} or shortly, \emph{deletion} and
\emph{insertion}, respectively.

\begin{dfn}{\rm
We say that a presentation is \emph{Peiffer aspherical} if every identity among relations can be carried to the empty
one by a sequence of Peiffer transformations.}
\end{dfn}

In particular, a presentation admitting identities among relations of odd length is automatically not Peiffer aspherical.

A large class of Peiffer aspherical presentations can be obtained by using Theorems~3.1 and 4.2, and Lemma~5.1
from~\cite{CCH}. They state, respectively, that Peiffer asphericity is preserved under certain $HNN$ extensions, under
free products, and under Tietze transformations.

In the next section we shall argue using Peiffer asphericity. However, for completeness, we mention that in the
literature there are (at least) three concepts of asphericity for presentations, which do not agree in general:
\emph{Peiffer asphericity} (called \emph{combinatorial asphericity} in~\cite{CCH}, see Proposition~1.5 there);
\emph{diagrammatical asphericity} defined in~\cite{CCH} like Peiffer asphericity but without allowing insertions (and
also considered in Chapter III.10 of~\cite{LS}); and \emph{topological asphericity}.

Let $H=\langle x_1,\ldots ,x_n\,|\, R_1,\ldots ,R_m \rangle$ be a presentation and $\mathcal{K}(H)$ be the
two-dimensional CW-complex with a single 0-cell, $n$ 1-cells corresponding to the generators $x_1,\ldots ,x_n$, and $m$
2-cells each one being attached to the 1-skeleton along the path determined by the spelling of the corresponding
relation. The presentation $H$ is said to be \emph{topologically aspherical} if $\pi_2(\mathcal{K}(H))=0$. As was
indicated in Proposition~1.1 of~\cite{CCH}, this is equivalent to the triviality of the second homology group of the
universal cover of $\mathcal{K}(H)$.

The relations between these three concepts are as follows (for more details, see the introduction and Proposition~1.3
of~\cite{CCH}):
\begin{itemize}
\item[(i)] topological asphericity implies Peiffer asphericity,
\item[(ii)] diagrammatical asphericity implies Peiffer asphericity,
\item[(iii)] for presentations where every relation is reduced, topological asphericity is equivalent to Peiffer
asphericity plus conciseness and ``no relator being a proper power".
\end{itemize}

\section{Proof of Theorem~\ref{pres}}\label{Mihailova}

Back to Mihailova's construction for $H=\langle x_1,\ldots ,x_n\,|\, R_1,\ldots ,R_m\rangle$, we recall that
$M(H)\leqslant F_n\times F_n$ is generated by $(x_i,x_i)$ and $(1,R_j)$, $i=1,\ldots ,n$, $j=1,\ldots ,m$. So, letting
$F_{n+m}$ be the free group with basis $\{ d_1,\ldots ,d_n, t_1,\ldots ,t_m \}$, we have an epimorphism $\pi \colon
F_{n+m}\to M(H)$ defined by $d_i \mapsto (x_i,x_i)$ and $t_j \mapsto (1,R_j)$, $i=1,\ldots ,n$, $j=1,\ldots ,m$. Now,
for proving Theorem~\ref{pres} we have to show that $\ker \pi$ is precisely the normal closure of the relations shown
in the pretended presentation for $M(H)$. Note that the images of elements $d_1,\dots ,d_n$ generate the diagonal
subgroup of $F_n\times F_n$, denoted $Diag(F_n\times F_n)$, which is isomorphic to $F_n$; hence, $\pi$ restricts to an
isomorphism from $D_n=\langle d_1,\ldots ,d_n\rangle \leqslant F_{n+m}$ onto $Diag(F_n\times F_n)\leqslant
M(H)\leqslant F_n\times F_n$.

We will keep the following notational convention in the proof: capital letters will always mean words on
$x_1,\ldots,x_n$; with this in mind, if $u$ is a word on $d_1,\ldots ,d_n$, then its capitalization $U$ will denote the
word obtained from $u$ by replacing each occurrence of $d_i$ to $x_i$. Thus, $U$ is just the projection of $\pi(u)$ to
the first (or the second) coordinate.

\medskip

{\it Proof of Theorem~\ref{pres}.} Recall that in the statement, $r_j$ is the word in $D_n$ obtained from $R_j$ by
replacing each $x_i$ to $d_i$, $j=1,\ldots ,m$.

Let $\mathcal{N}$ be the normal closure (in the free group $F_{n+m}$) of the recursive family of commutators
 $$
\{ [t_j, d^{-1}t_i^{-1}r_i \,d],\,\,\, [t_i,\,\, \text{\rm root}(r_i )] \mid i,j=1,\ldots ,m, \quad d\in D_n \}.
 $$
Our goal is to show that $\mathcal{N}=\ker \pi$. The inclusion $\mathcal{N}\leqslant \ker \pi$ is straightforward from
the following computations:
 $$
\pi([t_j, d^{-1}t_i^{-1}r_i\,d])=[(1,R_j), (u,u)^{-1}(R_i,1)(u,u)]=[(1,R_j), (u^{-1}R_iu,1)]=(1,1),
 $$
 $$
\pi([t_i,\,\, \text{\rm root}(r_i)]) =[(1,R_i), (\text{\rm root}(R_i),\text{\rm root}(R_i))]=(1,1).
 $$
In order to prove $\ker \pi \leqslant \mathcal{N}$, we shall use the following strategy: to each word $w\in \ker \pi$
we will associate an identity among relations for the presentation $\langle x_1,\dots ,x_n \mid R_1,\dots ,R_m\rangle$
of $H$, in such a way that if $w\neq 1$ then the associated identity is non-empty; then we will show that, after
applying an arbitrary Peiffer transformation, the resulting identity among relations is again the one associated to
some other word $w'\in \ker \pi$ satisfying, additionally, that $w^{-1}w' \in \mathcal{N}$.

Having seen this, let $w\in \ker \pi$ and consider the associated identity among relations. Since, by hypothesis, the
presentation $\langle x_1,\dots ,x_n \mid R_1,\dots ,R_m\rangle$ is Peiffer aspherical, there exists a sequence of
Peiffer transformations reducing such identity to the empty one. Now, repeatedly using the result mentioned in the
previous paragraph, we obtain a list of words (ending with the trivial one because the last identity is empty), $w,\,
w',\, w'',\ldots ,1$, and such that the difference between every two consecutive ones belongs to $\mathcal{N}$. This
shows that $w\in \mathcal{N}$ concluding the proof.

So, we are reduced to construct such an association. Let $w\in \ker \pi \leqslant F_{n+m}$ and write it in the form
$w=u_1t_{i_1}^{\varepsilon_1}u_2\cdots u_lt_{i_l}^{\varepsilon_l}u_{l+1}$, where $l\geqslant 0$ and $u_1,\ldots
,u_{l+1}$ are words in $d_1,\dots ,d_n$. Then, projecting $\pi(w)$ to each coordinate, we have
 \begin{equation}\label{1}
U_1U_2\cdots U_{l+1}=1 \quad \text{ and } \quad U_1R_{i_1}^{\varepsilon_1}U_2\cdots U_lR_{i_l}^{\varepsilon_l}
U_{l+1}=1.
 \end{equation}
Denote the accumulative products by $\mathbb{U}_i=U_1U_2\cdots U_i$, $i=1,\ldots ,l+1$ (note that
$\mathbb{U}_{l+1}=1$). By~(\ref{1}), we have
 $$
\mathbb{U}_1R_{i_1}^{\varepsilon_l}\mathbb{U}_1^{-1}\cdot \mathbb{U}_2R_{i_2}^{\varepsilon_2}\mathbb{U}_2^{-1}\cdot
\ldots \cdot \mathbb{U}_lR_{i_l}^{\varepsilon_l}\mathbb{U}_l^{-1}=1
 $$
in the free group $F_n$. In other words,
 \begin{equation}\label{2}
(\mathbb{U}_1R_{i_1}^{\varepsilon_l}\mathbb{U}_1^{-1}, \mathbb{U}_2R_{i_2}^{\varepsilon_2}\mathbb{U}_2^{-1},\ldots
,\mathbb{U}_lR_{i_l}^{\varepsilon_l}\mathbb{U}_l^{-1})
 \end{equation}
is an identity among relations for the presentation $\langle x_1,\ldots ,x_n\,|\, R_1,\ldots ,R_m\rangle$ of $H$. This
is the {\it identity associated to} $w\in \ker \pi$. Note that if this identity is empty, that is $l=0$, then $w=u_1\in
\langle d_1,\ldots ,d_n\rangle \cap \ker \pi$ and so $w=1$.

Let us analyze the situation when we apply an arbitrary Peiffer transformation to this identity.

\medskip

{\bf Case 1}: Consider the exchange which, for some $1\leqslant p\leqslant l-1$, replaces the consecutive terms
 $$
\mathbb{U}_p R_{i_p}^{\varepsilon_p}\mathbb{U}_p^{-1} \quad \text{ and } \quad \mathbb{U}_{p+1}
R_{i_{p+1}}^{\varepsilon_{p+1}} \mathbb{U}_{p+1}^{-1},
 $$
in (2), by the terms
 \begin{equation}\label{3}
\mathbb{U}_{p+1}R_{i_{p+1}}^{\varepsilon_{p+1}} \mathbb{U}_{p+1}^{-1} \quad \text{ and } \quad
(\mathbb{U}_{p+1}R_{i_{p+1}}^{-\varepsilon_{p+1}}\mathbb{U}_{p+1}^{-1} \mathbb{U}_p)
R_{i_p}^{\varepsilon_p}(\mathbb{U}_p^{-1} \mathbb{U}_{p+1}R_{i_{p+1}}^{\varepsilon_{p+1}}\mathbb{U}_{p+1}^{-1}),
 \end{equation}
respectively. We claim that the identity among relations obtained in this way is precisely the one corresponding to the
word
 $$
w'=v_1t_{i_1}^{\varepsilon_1}\cdots v_{p-1}t_{i_{p-1}}^{\varepsilon_{p-1}} v_p t_{i_{p+1}}^{\varepsilon_{p+1}} v_{p+1}
t_{i_p}^{\varepsilon_p}v_{p+2}t_{i_{p+2}}^{\varepsilon_{p+2}}\cdots v_lt_{i_l}^{\varepsilon_l}v_{l+1},
 $$
where
 $$
\begin{array}{lll} v_1=u_1, &  \hspace*{5mm} v_{p}=u_pu_{p+1}, &  \hspace*{5mm} v_{p+3}=u_{p+3},\\ \hspace*{5mm} \vdots
& \hspace*{5mm} v_{p+1}=r_{i_{p+1}}^{-\varepsilon_{p+1}}u_{p+1}^{-1}, & \hspace*{10mm} \vdots \\ v_{p-1}=u_{p-1},
\hspace*{5mm} & \hspace*{5mm} v_{p+2}=u_{p+1}r_{i_{p+1}}^{\varepsilon_{p+1}}u_{p+2}, & \hspace*{5mm} v_{l+1}=u_{l+1}.
\end{array}
 $$
And we also claim that $w^{-1}w' \in \mathcal{N}$. This second assertion is easy to verify since we can obtain back $w$
from $w'$ by permuting the two consecutive subwords $u_{p+1}t_{i_{p+1}}^{\varepsilon_{p+1}}
r_{i_{p+1}}^{-\varepsilon_{p+1}}u_{p+1}^{-1}$ and $t_{i_p}^{\varepsilon_p}$. But the commutator of these two words is
an element of $\mathcal{N}$: for $\varepsilon_{p+1}=-1$ this is immediate; and for $\varepsilon_{p+1}=1$ it follows
from the facts that, modulo $\mathcal{N}$, $t_{i_p}$ (and so $t_{i_p}^{\varepsilon_p}$) commutes with
$u_{p+1}(t_{i_{p+1}}^{-1} r_{i_{p+1}})^{\pm 1} u_{p+1}^{-1}$, but also $t_{i_{p+1}}$ commutes with $t_{i_{p+1}}^{-1}
r_{i_{p+1}}$ (and so, $t_{i_{p+1}}^{-1}$ with $r_{i_{p+1}}$). Therefore, $w'$ equals $w$ modulo $\mathcal{N}$.

To see the first part of the claim, let us capitalize the $v_i$'s:
 $$
\begin{array}{lll} V_1=U_1, &  \hspace*{5mm} V_{p}=U_pU_{p+1}, &  \hspace*{9mm} V_{p+3}=U_{p+3},\\ \hspace*{5mm}\vdots &
\hspace*{5mm} V_{p+1}=R_{i_{p+1}}^{-\varepsilon_{p+1}}U_{p+1}^{-1}, & \hspace*{14mm} \vdots \\ V_{p-1}=U_{p-1},
\hspace*{5mm} & \hspace*{5mm} V_{p+2}=U_{p+1}R_{i_{p+1}}^{\varepsilon_{p+1}}U_{p+2}, & \hspace*{9mm} V_{l+1}=U_{l+1}.
\end{array}
 $$
And let us compute the $\mathbb{V}_i=V_1V_2\cdots V_i$ 's:
 $$
\begin{array}{lll} \mathbb{V}_1=\mathbb{U}_1, &  \hspace*{5mm} \mathbb{V}_{p}=\mathbb{U}_{p+1}, &  \hspace*{9mm}
\mathbb{V}_{p+3}=\mathbb{U}_{p+3},\\ \hspace*{5mm}\vdots & \hspace*{5mm} \mathbb{V}_{p+1}=\mathbb{U}_{p+1}
R_{i_{p+1}}^{-\varepsilon_{p+1}}\mathbb{U}_{p+1}^{-1}\mathbb{U}_p, & \hspace*{14mm} \vdots \\ \mathbb{V}_{p-1}
=\mathbb{U}_{p-1},\hspace*{5mm} & \hspace*{5mm} \mathbb{V}_{p+2}=\mathbb{U}_{p+2}, & \hspace*{9mm}
\mathbb{V}_{l+1}=\mathbb{U}_{l+1}.
\end{array}
 $$
Finally, the identity among relations associated to $w'$ is
 $$
\begin{array}{rcl} (\mathbb{V}_1R_{i_1}^{\varepsilon_1}\mathbb{V}_1^{-1} & = & \mathbb{U}_1R_{i_1}^{\varepsilon_1}
\mathbb{U}_1^{-1},\\ & \vdots & \\ \mathbb{V}_{p-1}R_{i_{p-1}}^{\varepsilon_{p-1}}\mathbb{V}_{p-1}^{-1} & = &
\mathbb{U}_{p-1}R_{i_{p-1}}^{\varepsilon_{p-1}}\mathbb{U}_{p-1}^{-1},\\ \mathbb{V}_p R_{i_{p+1}}^{\varepsilon_{p+1}}
\mathbb{V}_p^{-1} & = & \mathbb{U}_{p+1}R_{i_{p+1}}^{\varepsilon_{p+1}}\mathbb{U}_{p+1}^{-1},\\ \mathbb{V}_{p+1}
R_{i_p}^{\varepsilon_p}\mathbb{V}_{p+1}^{-1} & = & \mathbb{U}_{p+1}R_{i_{p+1}}^{-\varepsilon_{p+1}}
\mathbb{U}_{p+1}^{-1} \mathbb{U}_pR_{i_p}^{\varepsilon_p} \mathbb{U}_p^{-1} \mathbb{U}_{p+1}
R_{i_{p+1}}^{\varepsilon_{p+1}} \mathbb{U}_{p+1}^{-1}, \\ \mathbb{V}_{p+2} R_{i_{p+2}}^{\varepsilon_{p+2}}
\mathbb{V}_{p+2}^{-1} & = & \mathbb{U}_{p+2}R_{i_{p+2}}^{\varepsilon_{p+2}}\mathbb{U}_{p+2}^{-1},\\ & \vdots & \\
\mathbb{V}_lR_{i_l}^{\varepsilon_l}\mathbb{V}_l^{-1} & = & \mathbb{U}_lR_{i_l}^{\varepsilon_l}\mathbb{U}_l^{-1}),
\end{array}
 $$
which does coincide with the identity among relations obtained from~(\ref{2}) after applying the Peiffer
transformation~(\ref{3}).

\medskip

{\bf Case 2}: Consider the deletion which, for some $1\leqslant p\leqslant l-1$, deletes the consecutive terms
 \begin{equation}\label{4}
\mathbb{U}_pR_{i_p}^{\varepsilon_p}\mathbb{U}_p^{-1} \quad \text{ and } \quad \mathbb{U}_{p+1}
R_{i_{p+1}}^{\varepsilon_{p+1}} \mathbb{U}_{p+1}^{-1},
 \end{equation}
in~(\ref{2}), assuming that its product equals 1. We claim that the identity among relations obtained in this way is
precisely the one corresponding to the word
 $$
w'=v_1t_{i_1}^{\varepsilon_1}\cdots v_{p-1}t_{i_{p-1}}^{\varepsilon_{p-1}} v_p t_{i_{p+2}}^{\varepsilon_{p+2}} v_{p+1}
\cdots v_{l-2} t_{i_l}^{\varepsilon_l} v_{l-1},
 $$
where
 $$
\begin{array}{lll} v_1=u_1, &  \hspace*{5mm} &  \hspace*{5mm} v_{p+1}=u_{p+3},\\ \hspace*{5mm}\vdots & \hspace*{5mm}
v_{p}=u_pu_{p+1}u_{p+2}, & \hspace*{10mm} \vdots \\ v_{p-1}=u_{p-1},\hspace*{5mm} & \hspace*{5mm} & \hspace*{5mm}
v_{l-1}=u_{l+1}.
\end{array}
 $$
And we also claim that $w^{-1}w' \in \mathcal{N}$. This second assertion follows from the hypothesis that
$(\mathbb{U}_p R_{i_p}^{\varepsilon_p} \mathbb{U}_p^{-1})\cdot (\mathbb{U}_{p+1} R_{i_{p+1}}^{\varepsilon_{p+1}}
\mathbb{U}_{p+1}^{-1})=1$. In fact, conciseness implies that $i_p=i_{p+1}$, $\varepsilon_p =-\varepsilon_{p+1}$ and so
$\mathbb{U}_p^{-1}\mathbb{U}_{p+1} =U_{p+1}$ commutes with $R_{i_{p+1}}$; hence, $u_{p+1}$ commutes with $r_{i_{p+1}}$
and so $u_{p+1}\in \langle \text{\rm root}(r_{i_{p+1}})\rangle$. Now $w'$ can be obtained from $w$ by replacing the
subword $t_{i_p}^{\varepsilon_p}u_{p+1} t_{i_{p+1}}^{\varepsilon_{p+1}}$ to $u_{p+1}$. But
$(t_{i_p}^{\varepsilon_p}u_{p+1} t_{i_{p+1}}^{\varepsilon_{p+1}})^{-1}u_{p+1}\in \mathcal{N}$
since $t_{i_{p+1}}$ commutes with $\text{\rm root}(r_{i_{p+1}})$ modulo
$\mathcal{N}$.

To see the first part of the claim, let us capitalize the $v_i$'s:
 $$
\begin{array}{lll} V_1=U_1, &  \hspace*{5mm} &  \hspace*{5mm} V_{p+1}=U_{p+3},\\ \hspace*{5mm}\vdots & \hspace*{5mm}
V_{p}=U_pU_{p+1}U_{p+2}, & \hspace*{10mm} \vdots \\ V_{p-1}=U_{p-1},\hspace*{5mm} & \hspace*{5mm} & \hspace*{5mm}
V_{l-1}=U_{l+1}.
\end{array}
 $$
And let us compute the $\mathbb{V}_i=V_1V_2\cdots V_i$ 's:
 $$
\begin{array}{lll} {\mathbb{V}}_1={\mathbb{U}}_1, &  \hspace*{5mm} &  \hspace*{5mm} {\mathbb{V}}_{p+1}={\mathbb{U}}_{p+3},
\\ \hspace*{5mm}\vdots & \hspace*{5mm} {\mathbb{V}}_{p}={\mathbb{U}}_{p+2}, & \hspace*{10mm} \vdots \\ {\mathbb{V}}_{p-1}
={\mathbb{U}}_{p-1},\hspace*{5mm} & \hspace*{5mm} & \hspace*{5mm} {\mathbb{V}}_{l-1}={\mathbb{U}}_{l+1}.
\end{array}
 $$
Finally, the identity among relations associated to $w'$ is
 $$
\begin{array}{rcl} (\mathbb{V}_1R_{i_1}^{\varepsilon_1}\mathbb{V}_1^{-1} & = & \mathbb{U}_1R_{i_1}^{\varepsilon_1}
\mathbb{U}_1^{-1},\\ & \vdots & \\ \mathbb{V}_{p-1}R_{i_{p-1}}^{\varepsilon_{p-1}}\mathbb{V}_{p-1}^{-1} & = &
\mathbb{U}_{p-1}R_{i_{p-1}}^{\varepsilon_{p-1}}\mathbb{U}_{p-1}^{-1},\\ \mathbb{V}_pR_{i_{p+2}}^{\varepsilon_{p+2}}
\mathbb{V}_p^{-1} & = & \mathbb{U}_{p+2}R_{i_{p+2}}^{\varepsilon_{p+2}}\mathbb{U}_{p+2}^{-1},\\ \mathbb{V}_{p+1}
R_{i_{p+3}}^{\varepsilon_{p+3}}\mathbb{V}_{p+1}^{-1} & = & \mathbb{U}_{p+3}R_{i_{p+3}}^{\varepsilon_{p+3}}
\mathbb{U}_{p+3}^{-1},\\ & \vdots & \\ \mathbb{V}_{l-2} R_{i_l}^{\varepsilon_l}\mathbb{V}_{l-2}^{-1} & = & \mathbb{U}_l
R_{i_l}^{\varepsilon_l}\mathbb{U}_l^{-1}), \end{array}
 $$
which coincides with the identity among relations obtained from~(\ref{2}) after applying the Peiffer
transformation~(\ref{4}).

\medskip

\textbf{Case 3:} Consider an insertion, and argue in a similar way as in Case~2.

\medskip

This concludes the proof. \qed

\section{A recursively presented orbit undecidable subgroup of $\aut(F_3)$}\label{OUD}

In~\cite{BMV}, O.~Bogopolski, A.~Martino and E.~Ventura studied the conjugacy problem for extensions of groups. In that
context, the notion of orbit decidability is crucial and we recall it here.

Let $F$ be a group, and $A\leqslant \aut(F)$. We say that $A$ is \emph{orbit decidable} if and only if there exists an
algorithm such that, given $u,v\in F$, decides whether $v$ is conjugate to $\alpha (u)$ for some $\alpha \in A$.

The main result in~\cite{BMV} states that, given a short exact sequence of groups
 $$
1\to F\to G\to P\to 1
 $$
with some conditions on $F$ and $P$, the group $G$ has solvable conjugacy problem if and only if the action subgroup
 $$
A_G=\{ \gamma_g \colon F\to F,\, x\mapsto g^{-1}xg \mid g\in G\}\leqslant \aut(F)
 $$
is orbit decidable (see~\cite[Theorem~3.1]{BMV} for details).

In particular, this applies to the case where $F$ and $P$ are finitely generated free groups, giving a characterization
of the solvability of the conjugacy problem within the family of [f.g. free]-by-[f.g. free] groups. This family of
groups is interesting because C.F.~Miller, back in the 1970's, already showed the existence of
[f.g.~free]-by-[f.g.~free] groups with unsolvable conjugacy problem (see~\cite{M1}). Via~\cite[Theorem~3.1]{BMV}, this
can be restated by saying that $\aut (F_n)$ contains finitely generated orbit undecidable subgroups (for some $n$).

Question~6 in the last section of~\cite{BMV} asks whether finitely presented subgroups $A\leqslant \aut(F_n)$ are orbit
decidable or not. The answer is known to be positive in rank 2 (every finitely generated subgroup of $\aut (F_2)$ is
orbit decidable, see~\cite[Proposition~6.21]{BMV}), but open for bigger rank. The comment made in~\cite{BMV} after this
question says that if $H$ is a finitely generated group with unsolvable word problem, then Mihailova's group $M(H)$ is
isomorphic to an orbit undecidable subgroup of $\aut(F_3)$. And, as mentioned in the introduction, this subgroup is
then finitely generated, and recursively presented, but it cannot be finitely presented.

In the rest of the paper, we will recall how $M(H)$ can be embedded into $\aut(F_3)$, in such a way that the image
becomes an orbit undecidable subgroup of $\aut(F_3)$. Then we will choose an appropriate $H$ and prove
Theorem~\ref{corol} by applying Theorem~\ref{pres} to $A=M(H)$ .

Of course, Theorem~\ref{corol} does not answer the above mentioned Question~6, but shows its tightness in the sense
that orbit undecidability is already showing up in the class of one-parametric recursively presented subgroups of
$\aut(F_3)$.

First, let $F_3 =\langle q,a,b \mid \quad\rangle$ be the free group on $\{q,a,b\}$, and let us embed $F_2 \times F_2$
into $\aut(F_3)$ in the following natural way. For every $u,v\in \langle a,b\rangle$, consider the automorphism
 $$
\begin{array}{rcl} _u\theta_v \colon F_3 & \to & F_3 \\ q & \mapsto & uqv \\ a & \mapsto & a \\ b & \mapsto & b.
\end{array}
 $$
Clearly, $_{u_1}\theta_1 \cdot _{u_2}\theta_1 =\, _{u_1u_2}\theta_1$ and $_1 \theta_{v_1} \cdot _1\theta_{v_2} =\,
_1\theta_{v_2v_1}$, which means that $\{ \, _u\theta_1 \mid u\in \langle a,b\rangle\} \simeq F_2$ and $\{ \, _1\theta_v
\mid v\in \langle a,b\rangle\} \simeq F_2^{\rm op}\simeq F_2$. It is also clear that $_u\theta_1 \cdot \,_1\theta_v =\,
_u\theta_v =\, _1\theta_v \cdot _u\theta_1$. So, we have an embedding $\theta \colon F_2\times F_2 \simeq F_2^{\rm
op}\times F_2^{\rm op} \hookrightarrow \aut(F_3)$ given by $(u, v)\mapsto \, _{u^{-1}}\theta_v$, whose image is
 $$
F_2 \times F_2 \simeq B=\langle _{a^{-1}}\theta_1 ,\, _{b^{-1}}\theta_1,\, _1\theta_a ,\, _1\theta_b \rangle =\{ \,
_u\theta_v \mid u,\, v\in \langle a,b\rangle\} \leqslant \aut(F_3).
 $$
As shown in~\cite[Section~7.2]{BMV}, the element $qaqbq$ satisfies the technical condition required
in~\cite[Proposition~7.3]{BMV}. Hence, we have

\begin{lem}[7.3 in~\cite{BMV}]\label{a-b-aut}
For the above defined subgroup $B\leqslant \aut(F_3)$ and for every subgroup $A\leqslant B$, undecidability of the
membership problem for $A$ in $B$ implies orbit undecidability for $A$ in $\aut(F_3)$.
\end{lem}

We are ready to deduce the main result of this section.

\begin{thm}\label{corol}
There exists a finitely generated (and not finitely presented) orbit undecidable subgroup $A\leqslant \aut(F_3)$
admitting a one-parametric recursive presentation as in Theorem~\ref{pres}.
\end{thm}

\demo In \cite{CM1}, D.J.~Collins and C.F.~Miler~III proved that there exists a finite, concise and Peiffer aspherical
presentation $\langle x_1,\dots ,x_n \,|\, R_1,\dots ,R_m\rangle$ of a group $H$ with unsolvable word problem. The
corresponding Mihailova's group $M(H)$ is a subgroup of $F_n\times F_n$ and the membership problem for $M(H)$ in
$F_n\times F_n$ is unsolvable.

Now, denoting $A=M(H)$ and using a finite index embedding of $F_n\times F_n$ in $B\cong F_2\times F_2$, we have that
$A\leqslant B$ and the membership problem for $A$ in $B$ is unsolvable. By Lemma~\ref{a-b-aut}, $A$ is an orbit
undecidable subgroup of $\aut(F_3)$.

Moreover, as it was discussed in the introduction, $A$ is finitely generated, and is not finitely presented. But
Theorem~\ref{pres} provides an explicit one-parametric recursive presentation for~$A$. This concludes the proof. \qed

\medskip

We end by reproducing \cite[Question~6]{BMV} again:

\begin{qs}
Does there exist a finitely presented orbit undecidable subgroup of $\aut(F_n)$, for~$n\geqslant~\!3$\,?
\end{qs}


\section{Acknowledgements}

The first named author thanks the MPIM at Bonn for its support and excellent working conditions during the fall 2008,
while this research was finished. The second author gratefully acknowledges partial support from the MEC (Spain) and
the EFRD (EC) through project number MTM2006-13544.

\end{document}